 \documentclass[12pt]{article}
 \usepackage{graphicx}
 \usepackage[english]{babel}
 \usepackage{amsmath,amssymb,euscript,amsthm,amsfonts,mathrsfs,amscd}
 \usepackage{color}
 \usepackage{floatflt}
 \usepackage{mathrsfs} 

 \theoremstyle{plain}

 \newtheorem{prop}{Proposition}
 \newtheorem{Def}{Definition}
 \theoremstyle{definition}
 \newtheorem{st}{Step}

 \theoremstyle{remark}
  \newtheorem{rem}{Remark}

\newcommand{\dd}{\ensuremath{\displaystyle}}

\begin{document}
{\selectlanguage{english}
\title{Lorden's inequality and coupling method for backward renewal process}

\author{G. A. Zverkina}
\sloppy

\maketitle
\begin{abstract}
We give a scheme of  using  the coupling method to obtain strong bounds for the convergence rate of the distribution of the backward renewal process in the total variation distance.
This scheme can be applied to a wide class of regenerative processes in queuing theory.
\\
\textbf{keywords ~ }{\emph{backward renewal process, renewal process, convergence rate, strong bounds, total variation metric, Lorden's inequality.}}
\end{abstract}

\section{Introduction}
\label{sec:intro}
Obviously, the behaviour of some queueing system (or of some reliability system) can be described by regenerative process.
Hence a study of the behaviour of regenerative processes is an important problem in the queuing theory and in the reliability theory.

\begin{Def} Recall that the stochastic process $\{X_t,\,t\geqslant0\}$ with the measurable state space $(\mathcal{X},\mathcal{B}(\mathcal{X}))$ and filtration $\{\mathcal{F}_t,\,t\geqslant 0\}$ defined on the probability space $(\Omega, \mathcal F, \mathbf{P})$ is called regenerative process if there exists a sequence of  Markov moments $\{\theta_i\}_{i\in \mathbb N}$  with respect to the filtration $\mathcal{F}_t$ such that:

  1. $X_{\theta_i}=X_{\theta_j}$ for all $i,j\in\mathbb N$;

  2. The random elements $\Xi_i \stackrel{{\rm def}}{=\!\!\!=} \{X_t, \, t\in[\theta_i,\theta_{i+1}]\}$ ($i\in\mathbb N$), are mutually independent and identically distributed. \hfill\ensuremath{\rhd}
\end{Def}
Hence, the random variables $\zeta_{i+1} \stackrel{{\rm def}}{=\!\!\!=} \theta_{i+1}-\theta_i$, $i\in\mathbb N$ are i.i.d.; denote

$$F(s) \stackrel{{\rm def}}{=\!\!\!=}  \mathbf{P}\{\zeta_i\leqslant s\},\;\; i\in\mathbb N,$$ and 
$$G(s) \stackrel{{\rm def}}{=\!\!\!=}  \mathbf{P}\{\theta_1\leqslant s\}.$$

Also, the random variables $\big\{\{\zeta_i\}_{i=1}^\infty$ and $\theta_1\big\}$ are mutually independent.

If the random variables $\zeta_i$, $i\in\mathbb N$, and $\zeta_1 \stackrel{{\rm def}}{=\!\!\!=} \theta_1$ have finite expectations, then the ergodicity of the process $X_t$ follows from the Harris-Khasminsky principle, i.e. the distribution of $X_t$ weekly converges to the unique probability ({\it stationary}) measure $\mathcal{P}$ on $(\mathcal{X},\mathcal{B}(\mathcal{X}))$.

Moreover, in the 70s A.A.~Borovkov \cite{borovkov} showed that if for some $\kappa>1$ the condition
$$  \big\{\mathbf E\,\theta_1^\kappa<\infty, \;\; \mathbf{E}\,\zeta_2^\kappa<\infty\big\}$$
is satisfied, then for all $\alpha\leqslant \kappa-1$ there exists the (unknown) constant $K(\alpha)$ such that for all $S\in\mathcal{B}(\mathcal{X})$ and all $t>0$ the inequality
$$|\mathbf{P}\{X_t\in S\}-\mathcal{P}(S) |<K(\alpha)t^{-\alpha}$$ 
is true.

The bounds for the constant $K(\alpha)$ was defined in some particular cases -- see e.g.\cite{zvbib13,zvbib14,zvbib15,zvbib16} et al.; our goal is to give some procedure to obtain such bounds for a sufficiently wide class of regenerative processes.

Firstly we remark that the Markov moments  $\{\theta_i\}_{i\in \mathbb N}$ form an embedded renewal process $R_t \stackrel{{\rm def}}{=\!\!\!=} \sum\limits_{i=1}^\infty \mathbf{1}(\theta_i\leqslant t)$.

If $\theta_1=0$, then $R_t$ is a renewal process without delay.

If $\theta_1>0$, then we can interpret the time interval $(0,\theta_1)$ as the residual time of the regeneration period of $X_t$ (or the renewal period of $R_t$), where the process $X_t$ began its regeneration period at the some fixed time $(-a)$.

The random variable $D_t \stackrel{{\rm def}}{=\!\!\!=}  (\theta_{R_t+1}-t)$ is called {\it forward renewal time}.
For $t>\theta_1$,   \emph{Lorden's inequality} \cite{lorden} is true:
\begin{equation}\label{lorden}
  \mathbf E\,D_t\leqslant \Theta \stackrel{{\rm def}}{=\!\!\!=}  \frac{\mathbf E\,\zeta_2^2}{\mathbf E\,\zeta_2}=\frac{\dd\int\limits _0^\infty s^2 \,\mathrm{d} \, F(s)}{\dd\int\limits _0^\infty (1-F(s)) \,\mathrm{d} \, s}.
\end{equation}

Also denote $B_t \stackrel{{\rm def}}{=\!\!\!=}  t-\theta_{R_t}$ for $t\geqslant\theta_1$; $B_t$ is a time from the last renewal ($\theta_{R_t}\leqslant t$) to the time $t$, it is called {\it backward renewal time} of the renewal process $R_t$.

The state space of the process $B_t$ is $(\mathbb R_+,\mathcal{B}(\mathbb R_+))$; if $\theta_1=0$ then the process  $B_t$ starts from the state $B_0=0$ (non-delay process).

If $\theta_1\neq 0$ then we consider the random variable $\theta_1$ as the residual time of the first renewal period of the renewal process $R_t$ started at some fixed time $(-a)$; consequently we put $$G(s)=F_a(s) \stackrel{{\rm def}}{=\!\!\!=}  \displaystyle \frac{F(a+s)-F(a)}{1-F(a)};$$ 
{\it we suppose that $F(a)<1$ for all $a\in \mathbb R$}.
So, for $t\in[0,\theta_1]$ we put $B_t=a+t$, in the assumption that $G(s)=F_a(s)$.

It is easy to see that the process $B_t$ is Markov.
If $\mathbf E\,\zeta_2<\infty$, and consequently $\mathbf E\,\theta_1<\infty$, then $B_t$ is ergodic, i.e. its distribution weekly converges to the stationary distribution as $t\to \infty$.

Moreover, if the distribution of $B_t$ converges to the stationary distribution, then the distribution of  $X_t$ also converges to the stationary distribution as  $t\to \infty$. 

And if we know the bounds of the convergence rate of distribution of the process $B_t$ (in some sense), then we know the bounds of the convergence rate of distribution of the process $X_t$ (in the same sense), because the distribution of $X_t$ is determined by the value of $B_t$.

So, our goal is to obtain the bounds for the convergence rate of the backward renewal process, and for this aim we will use the {\it coupling method}.

\section{Coupling method}
\label{sec:base-section}
Let $X_t'$ and $X_t''$ be the homogeneous independent Markov processes with  the same state space $(\mathcal{X},\mathcal{B}(\mathcal{X}))$ and the same transition function, but with different initial states: $X_0'=x'\neq X_0''=x''$.
Denote the distribution of the process $X_t$ with the initial state $x$ at the time $t$ by $\mathcal{P}_t^x$, i.e.
$$\mathcal{P}_t^{x}(S)=\mathbf{P}\{X_t\in S|X_0=x\};$$
and for all $x\in\mathcal{X}$, $\mathcal{P}_t^x\Longrightarrow \mathcal{P}$ as $t\to \infty$.

Let the paired process $\mathcal{Z}_t=\left(U_t',U_t''\right)$ created on some probability space satisfies the following conditions $(i)$--$(iii)$:

$(i)$ For all $t\geqslant 0$ and $S\in\mathcal{B}(\mathcal{X})$, $\mathbf{P}\left\{U_t'\in S\right\}= \mathbf{P}\left\{X'_t\in S\right\}$, and $\mathbf{P}\left\{U_t''\in S\right\}= \mathbf{P}\left\{X''_t\in S\right\}$,  therefore, $U'_0=X'_0=x'$ and $U''_0=X''_0=x''$.

$(ii)$ For all $t\geqslant \tau\left(x',x''\right) \stackrel{{\rm def}}{=\!\!\!=}  \inf\left\{t\geqslant 0: \, U'_t=U''_t\right\}$, the equality $U'_t=U''_t$ is true.

$(iii)$ For all $x'$, $x''\in \mathcal{X}$, $\mathbf{P}\{\tau(x','')<\infty\}=1$.

The paired process $\mathcal{Z}_t$ satisfying conditions  $(i)$--$(iii)$  is called {\it successful coupling} of the processes $X_t'$ and $X_t''$, and for them, the based coupling inequality can be written so: for all $S\in\mathcal{B}(\mathcal{X})$ 
\begin{multline*}
  \left|\mathcal{P}_t^{x'}(S)-\mathcal{P}_t^{x'' }(S)\right|=
  \\=\big|\mathbf{P}\{X_t' \in S\}-\mathbf{P}\{X_t '' \in S\}\big|=\big|\mathbf{P}\{U_t' \in S\}-\mathbf{P}\{U_t''\in S\}\big|=
\\ \\
=\big|\mathbf{P}\{U_t' \in S\,\&\,\tau(x',x'')\leqslant t\}+\mathbf{P}\{U_t' \in S\,\&\,\tau(x',x'')> t\}-
\\ \\
-(\mathbf{P}\{U''_t\in S\,\&\,\tau(x',x'')\leqslant t\} +\mathbf{P}\{U''_t\in S\,\&\,\tau(x',x'')> t\})\big|\leqslant
\\ \\
  \leqslant \big|\mathbf{P}\{U_t' \in S\,\&\,\tau(x',x'')> t\}-\mathbf{P}\{U''_t \in S\,\&\,\tau(x',x'')> t\}\big|\leqslant
\\ \\
\leqslant\mathbf{P}\{\tau(x',x'')> t\},
\end{multline*}
and if we can find an estimate $\mathbf{P}\{\tau(x,y)> t\} \leqslant \varphi(x,y,t)$, and 
$$\widehat\varphi(x,t) \stackrel{{\rm def}}{=\!\!\!=}  \mathbf E\,\varphi\left(x,\widetilde X,t\right)<\infty
$$
for the random variable $\widetilde X$ with the stationary distribution $\mathcal{P}$, then for all $S\in\mathcal B(\mathcal X)$
\begin{equation}\label{oc}
  \left|\mathcal{P}_t^{x}(S)-\mathcal{P}(S)\right| \leqslant\int\limits _{\mathcal{X}}\varphi(x,u,t)\mathcal{P}(\,\mathrm{d}\, u)=\widehat\varphi(x,t),
\end{equation}
and
\begin{equation}\label{oc1}
\|\mathcal{P}_t^{x }-\mathcal{P}\|_{TV} \stackrel{{\rm def}}{=\!\!\!=}  2 \sup\limits_{S\in\mathcal{B}(\mathcal{X})} \left|\mathcal{P}_t^{X_0 }(S)-\mathcal{P}(S)\right|\leqslant 2\widehat\varphi(x,t).
\end{equation}
\section{Auxiliary considerations}

\begin{Def} \label{def2}The common part of the distributions of the random variables $\xi_1$ and $\xi_2$ with distribution functions $\Psi_j(s)=\mathbf{P}\{\xi_j\leqslant s\}$, $s\in \mathbb{R}$ (here and hereafter $j=1,2$) is

$$\varkappa \stackrel{{\rm def}}{=\!\!\!=} \varkappa(\Psi_1(s),\Psi_2(s)) \stackrel{{\rm def}}{=\!\!\!=} \displaystyle \int\limits _{-\infty}^\infty \min(\psi_1(u),\psi_2(u))\,\mathrm{d}\, u,$$ 
where $\psi_j(s) \stackrel{{\rm def}}{=\!\!\!=}  \Psi_j'(s)$ if $\Psi_j'(s)$ exists, and $\psi_j(s) \stackrel{{\rm def}}{=\!\!\!=}  0$ otherwise.\hfill\ensuremath{\rhd}
\end{Def}

\begin{prop}\label{prop1}
~

1.
If $\varkappa>0$, then the function
$
\widetilde\Psi(s) \stackrel{{\rm def}}{=\!\!\!=} \displaystyle \frac{1}{\varkappa}\displaystyle \int\limits _{-\infty}^s \min(\psi_1(u),\psi_2(u))\,\mathrm{d}\, u
$
is a distribution function, and the functions $\widetilde\Psi_j(s) \stackrel{{\rm def}}{=\!\!\!=} \displaystyle \frac{\Psi_j(s)-\varkappa\widetilde\Psi(s)}{1-\varkappa}$ are a distribution functions ($j=1,2$).

2.
Let $\mathcal{U}'$, $\mathcal{U}''$ and $\mathcal{U}'''$ be independent random variables with uniform distribution on $[0,1]$.
Then
$$
\widetilde\xi_j \stackrel{{\rm def}}{=\!\!\!=}  \widetilde\Psi^{-1}(\mathcal{U}'')\mathbf{1}(\mathcal{U}'<\varkappa) +\widetilde\Psi_j^{-1}(\mathcal{U}''')\mathbf{1}(\mathcal{U}'\geqslant\varkappa) \stackrel{\mathcal{D}}{=}  \xi_j,
$$
and  $$\mathbf{P}\left\{\widetilde\xi_1=\widetilde\xi_2\right\}=\varkappa.$$ 

{\it Here  and  hereafter  where  we  put  $h^{-1}(s) \stackrel{{\rm def}}{=\!\!\!=}  \inf\{x\in\mathbb{R}:\,h(x)\geqslant s\}$ for non-decreasing function $h(x)$}.\hfill\ensuremath{\rhd}
\end{prop}
\begin{rem}
Proposition \ref{prop1} is a simplified variant of the {\it Coupling Lemma} or ``{\it Lemma about three random variables}'' (see, e.g., \cite{ver}).\hfill\ensuremath{\rhd}
\end{rem}

\section{Successful coupling for backward renewal process}
Now we consider the backward renewal process $B_t$; its state space is $(\mathbb{R}_+,\mathcal{B}(\mathbb{R}_+))$; its stationary distribution has a distribution function

$$\widetilde{F}(s)=(\mathbf E\,\zeta_1)^{-1}\displaystyle \int\limits _0^s(1-F(u))\,\mathrm{d}\, u $$
and the stationary distribution is
$$
\mathcal P(S)=(\mathbf E\,\zeta_1)^{-1} \displaystyle \int\limits _{S} (1-F(u))\,\mathrm{d}\, u .
$$

\subsection{Basic assumption.}  Here and hereafter we suppose: $\kappa\geqslant 2$, and for some $A>0$ for all $t>A$ and for all $\varepsilon>0$,
$$
 \int\limits _t^{t+\varepsilon} f(s)\,\mathrm{d}\, s>0, \mbox{ where } f(s) \stackrel{{\rm def}}{=\!\!\!=} \begin{cases}
F'(s), & \mbox{if } \exists \,F'(s);
\\
0, & \mbox{otherwise},
\end{cases}
\mbox{ and }\mathbf E\,\zeta_2^\kappa<\infty.
\eqno(\ast)
$$

We will construct the successful coupling for two versions of the process $B_t$ started from different initial states $b_1$ and $b_2$; denote them $B_t^{(1)}$ and $B_t^{(2)}$.

\subsection{Construction of renewal process.}\label{constr}
Recall the method of construction of the renewal process $R_t$ with delay $\theta_1$ having the distribution function $G(s)$.

Let $\{\mathcal{U}_n\}$ be a sequence of independent random variables with uniform distribution on $[0,1]$.

The construction of the renewal times for the renewal process $R_t$, i.e. $\theta_1$, $\theta_2$, \ldots, $\theta_n$, \ldots \, is follow:\\
$\theta_1 \stackrel{{\rm def}}{=\!\!\!=}  G^{-1}(\mathcal{U}_1)=F_a^{-1}(\mathcal{U}_1)=\zeta_1;\; \theta_2 \stackrel{{\rm def}}{=\!\!\!=}  \theta_1^a+F^{-1}(\mathcal{U}_2)=\theta_1+\zeta_2;\ldots$

\hspace{4cm}$\ldots\, \theta_n \stackrel{{\rm def}}{=\!\!\!=}  \theta_{n-1}+F^{-1}(\mathcal{U}_n)=\theta_{n-1}+\zeta_n;\,\ldots$

So, for construction of two independent backward renewal processes we need two sequences of independent random variables with uniform distribution on $[0,1]$ -- let they be  $\{\mathcal{U}_{n,1}\}$ and  $\{\mathcal{U}_{n,2}\}$.
Now we denote for $j=1,2$
\begin{multline*}
 \theta_1^{(j)} \stackrel{{\rm def}}{=\!\!\!=}  F_{b_j}^{-1}(\mathcal{U}_{1,j})=\zeta_1^{(j)};\; \theta_2^{(j)} \stackrel{{\rm def}}{=\!\!\!=}  \theta_1^{(j)}+F^{-1}(\mathcal{U}_{2,j})=\theta_1^{(j)}+ \zeta_2^{(j)};\ldots
 \\ \\
 \ldots\, \theta_n^{(j)} \stackrel{{\rm def}}{=\!\!\!=}  \theta_{n-1}^{(j)}+F^{-1}(\mathcal{U}_{n,j})=\theta_{n-1}^{(j)}+ +\zeta_n^{(j)};\,\ldots,
\end{multline*}

\begin{equation}\label{eq:bt}
B^{(j)}_t=\begin{cases}
b_j+t, & \mbox{ if } t<\theta_1^{(j)};
\\ \\
t-\sum\limits_{i=1}^\infty\left( \mathbf{1}\left(\theta_i^{(j)}\leqslant t\right)\zeta_i^{(j)}\right),&\mbox{ otherwise};
\end{cases}
\end{equation}
and
\begin{equation}\label{eq:bt2}
R_t^{(j)}=\displaystyle \sum\limits_{i=1}^\infty\mathbf{1} \left(\theta_i^{(j)}\leqslant t\right ); \qquad D_t^{(j)}=\theta^{(j)}_{R_t^{(j)}+1}-t.
\end{equation}

From (\ref{eq:bt}) and (\ref{eq:bt2}) we see that the processes $B_t^{(1)}$ and $B_t^{(2)}$ are piecewise linear, and they can begin to be equal only at the time when both of them are equal to zero.
But for independent processes $B_t^{(1)}$ and $B_t^{(2)}$ the probability of they coincidence is zero because the distribution of residual time of any renewal periods of corresponding renewal process has a continuous component -- see ($\ast$).

Therefore we will construct the successful coupling concerning of two {\it dependent} processes $\widetilde{B}_t^{(1)}$ and $\widetilde{B}_t^{(2)}$ such that $\widetilde{B}_t^{(j)} \stackrel{\mathcal{D}}{=}  B_t^{(j)}$; for this aim we need an additional sequences of independent random variables $\{\mathcal{U}_{n}^i\}$ with uniform distribution on $[0,1]$.
\subsection{Construction of  the successful coupling \\ $\mathcal{Z}_t=\left(\widetilde{B}_t^{(1)}, \widetilde{B}_t^{(2)}\right)$.}\label{coup}
We will construct the paired process $\mathcal{Z}_t$ by following algorithm.

\begin{st}\label{st1} We begin to construct independent processes $\widetilde{B}_t^{(1)}$ and $\widetilde{B}_t^{(2)}$ according to the scheme \ref{constr} -- i.e. we construct the renewal times $\theta_i^{(j)}$ for the processes $\widetilde{B}_t^{(j)}$.
Put $T_1 \stackrel{{\rm def}}{=\!\!\!=}  \max\left(\theta_1^{(1)},\theta_1^{(2)}\right)$.

Note, that $T_1\leqslant\theta_1^{(1)}+\theta_1^{(2)}$.
For simplicity, here we suppose that $T_1=\theta_1^{(1)}$.
If $\theta_1^{(1)}=\theta_1^{(2)}$ then we go to the Step \ref{st3} $\Big($but  $ \mathbf{P}\left\{\theta_1^{(1)}=\theta_1^{(2)}\right\}=0\Big)$.
Otherwise, we go to the Step \ref{st2}.
\end{st}
\begin{st}\label{st2}
At the time $T_k=\theta_k^{(1)}\in\left(\theta^{(2)}_{\nu_k}, \theta^{(2)}_{\nu_k+1}\right)$, $\mathbf E\,D_{T_k}^{(2)}\leqslant \Theta$ (see (\ref{lorden})), and $\widetilde{B}_{T_k}^{(1)}=0\neq \widetilde{B}_{T_k}^{(2)}$.

By Markov inequality, for some $R>\Theta$, $$\mathbf{P}\left\{D_{T_k}^{(2)}\leqslant R\right\} \geqslant \pi_R \stackrel{{\rm def}}{=\!\!\!=} 1-\displaystyle \frac{\Theta}{R}.$$
And 
$$\mathbf{P} \left\{\theta_{k+1}^{(1)}-T_k =\zeta_{k}>R\right\}=1-F(R).$$

So, at the time $\theta_{\nu_k+1}^{(2)}=T_k+D_{T_k}^{(2)}$ with probability $P_R \stackrel{{\rm def}}{=\!\!\!=}  \pi_R(1-F(R))$ we have: $\widetilde{B}_{\theta_{\nu_{k}+1}^{(2)}}^{(2)}=0$ and $\beta \stackrel{{\rm def}}{=\!\!\!=}  \widetilde{B}_{\theta^{(2)}_{\nu_k+1}}^{(1)}=D_{T_k}^{(2)}\leqslant R$.

If $D_{T_k}^{(2)}> R$ or $\zeta_k^{(1)}<R$ then we move on to the next time $\theta_{k+1}^{(1)}$, i.e. we replace $k$ by $k+1$ and return to the Step \ref{st2}.

If $D_{T_k}^{(2)}\leqslant R$ and $\zeta_k^{(1)}\geqslant R$ then we stop both processes at the time $\theta_{\nu_k+1}^{(2)}$.
Then we prolong these processes (i.e. their residual times with distributions $F_\beta(s)$ and $F(s)$) using the additional random variables $\mathcal{U}_k^1$, $\mathcal{U}_k^2$, $\mathcal{U}_k^3$ --- by such a way that with probability $\varkappa(F_\beta (s), F(s))$ the next renewal times of both processes coincide (see  Proposition \ref{prop1}).

Note that the condition $(\ast)$ implies that the common part of distribution of the forward renewal times of both processes $$\varkappa(F_\beta(s),F(s))\geqslant \varkappa_R \stackrel{{\rm def}}{=\!\!\!=}  \inf\limits_{a\in[0,R]}\{\varkappa(F_a(s),F(s))\}>0$$ -- see Definition \ref{def2}.

Hence, at the time $\theta_{k+1}^{(1)}$ the constructed processes $\widetilde{B}_t^{(1)}$ and $\widetilde{B}_t^{(2)}$ coincide with probability $p_k\geqslant \varkappa_RP_R$ -- denote this event by $\mathcal{E}_k$; $ \mathbf{P}(\mathcal{E}_k)\geqslant \varkappa_RP_R$.

If $\widetilde{B}_{\theta_{k+1}^{(1)}}^{(1)}= \widetilde{B}_{\theta_{k+1}^{(1)}}^{(2)}(=0)$ the we go to the Step \ref{st3}.
Otherwise, we move on to the next time $\theta_{k+1}^{(1)}$, i.e. we replace $k$ by $k+1$ and return to the Step \ref{st2}.
\end{st}
\begin{st}\label{st3}
After the time $\tau \stackrel{{\rm def}}{=\!\!\!=} \theta_{k+1}^{(1)}$ such that $\widetilde{B}_{\tau}^{(1)}= \widetilde{B}_{\tau}^{(2)}=0$ we prolong the construction of the processes $\widetilde{B}_t^{(1)}$ and $\widetilde{B}_t^{(2)}$ identically by the scheme \ref{constr}.
Denote the event $\left\{\widetilde{B}_{\theta_{k+1}^{(1)}}^{(1)}= \widetilde{B}_{\theta_{k+1}^{(1)}}^{(2)}\,\& \, \widetilde{B}_{\theta_{k}^{(1)}}^{(1)}\neq \widetilde{B}_{\theta_{k}^{(1)}}^{(2)}\right\}$ by $\mathfrak{E}_k$; $$\mathfrak{E}_k=\left(\bigcap\limits_{m=2}^k\overline{\mathcal{E}}_m\right)\cap \mathcal{E}_{k+1};$$ 
$$ \mathbf{P}(\mathfrak{E}_k)=p_k\prod\limits_{m=2}^{k}(1-p_m)\leqslant (1-\varkappa_R P_R)^{k-1} \stackrel{{\rm def}}{=\!\!\!=}  q_R^{k-1}.
$$
\end{st}
\begin{prop}\label{prop2} The pared process $\mathcal{Z}_t=\left(\widetilde{B}_t^{(1)}, \widetilde{B}_t^{(2)}\right)$ is a successful coupling of processes ${B}_t^{(1)}$ and ${B}_t^{(2)}$. \hfill\ensuremath{\rhd}
\end{prop}
\section{Strong bounds for convergence rate.}
Now, for $\phi(s)=s^\alpha$ we can find an upper bound for 
$$ \mathbf E \,\phi(\tau(b_1,b_2))= \mathbf E \,(\tau(b_1,b_2)^\alpha), \;\; \alpha\geqslant 1$$
using  Jensen's  inequality  in  the  form 
$$\left(\sum\limits_{k=1}^n a_i\right)^\alpha \leqslant n^{\alpha-1}\left(\sum\limits_{k=1}^n a_i^\alpha\right)$$ 
for positive $a_i$, namely:
\begin{multline*}
    \mathbf E \,(\tau(b_1,b_2))^\alpha \leqslant  \sum\limits_{n=1}^\infty  \mathbf E \,\left(\left(T_1+ \sum\limits_{m=2}^{n+1} \zeta_m^{(1)}\right) \mathbf{1} (\mathfrak{E}_n)\right)^\alpha\leqslant 
    \\ \\
   \leqslant \sum\limits_{n=1}^\infty(n+2)^{\alpha-1} \mathbf E \,\left(\left(\left(\theta_1^{(1)}\right)^\alpha+ \left(\theta_1^{(2)}\right)^\alpha +  \sum\limits_{m=2}^{n+1} \left(\zeta_m^{(1)}\right)^\alpha\right) \mathbf{1} (\mathfrak{E}_n)\right)\leqslant 
   \\ \\
   \leqslant K_1(\alpha)\left( \mathbf E \,\left(\theta_1^{(1)}\right)^\alpha+ \mathbf E \,\left(\theta_1^{(2)}\right)^\alpha\right)+
   \\ \\
   + \sum\limits_{n=1}^\infty (n+2)^{\alpha-1}\left( \sum\limits_{m=2}^{n} \mathbf E \,\left(\left(\zeta_m^{(1)} \right)^\alpha \mathbf{1} (\overline{\mathcal{E}}_m)\right)q_R^{n-2}+\right.
   \\ \\
   + \mathbf E \,\left(\left(\zeta_{n+1}^{(1)}\right)^\alpha \mathbf{1} ({\mathcal{E}}_{n+1})\right)q_R^{n-1}\Bigg)\leqslant
   \\ \\
   \leqslant  K_1(\alpha)\left( \mathbf E \,\left(\theta_1^{(1)}\right)^\alpha+ \mathbf E \,\left(\theta_1^{(2)}\right)^\alpha\right)+ K_2(\alpha) \mathbf E \,\zeta_2^\alpha=\varpi_\alpha(b_1,b_2,t),
\end{multline*}
where  
$$K_1(\alpha)= \sum\limits_{n=1}^\infty (n+2)^{\alpha-1}q_R^{n-1},\;\; K_2(\alpha)= \sum\limits_{n=1}^\infty (n+2)^{\alpha}q_R^{n-1}, \;\; \sum\limits_{m=2}^1 \stackrel{{\rm def}}{=\!\!\!=}  0.$$

Now,  by  Markov  inequality,

$$ \mathbf{P}\{\tau(b_1,b_2)>t\}\leqslant\displaystyle  \frac{ \mathbf E \,(\tau(b_1,b_2))^\alpha}{t^\alpha}\leqslant \varphi_\alpha(b_1,b_2,t) \stackrel{{\rm def}}{=\!\!\!=} \frac{ \varpi_\alpha(b_1,b_2,t) }{t^\alpha} .
$$

Then, for $\alpha\leqslant \kappa-1$ we have
\begin{multline}\label{k}
\int\limits _{0}^\infty  \mathbf E \,(\tau(b_1,b_2))^\alpha  \,\mathrm{d} \, \widetilde F (b_2)\leqslant \int\limits _{0}^\infty \varpi_\alpha(b_1,b_2,t) \,\mathrm{d} \, \widetilde F (b_2)=
\\ \\
=  K_1(\alpha) \mathbf E \,\left(\theta_1^{(1)}\right)^\alpha+ K_2(\alpha) \mathbf E \,\zeta^\alpha +\int\limits _0^\infty K_1(\alpha)  \mathbf E \, \left(\theta_1^{(2)}\right)^\alpha \,\mathrm{d} \,\widetilde{F}(b_2)=
\\ \\
=K_1(\alpha) \mathbf E \,\left(\theta_1^{(1)}\right)^\alpha+ K_2(\alpha) \mathbf E \,\zeta^\alpha +\frac{K_1(\alpha)}{\alpha+1} \mathbf E \,\left(\theta_1^{(2)}\right)^{\alpha+1}= {K}(\alpha, b_1);
\end{multline}
and $\widehat{\varphi}_\alpha(b_1,t)= \displaystyle \frac{{K}(\alpha,b_1)}{t^\alpha}$ for $\alpha\in[1,\kappa-1]$.

Therefore we have the bounds (\ref{oc}) and (\ref{oc1}) for backward renewal process.
The estimate (\ref{k}) can be improved by the choice of $R$.
It is not optimal, and it can be done better by use of the properties of the distribution $F$, and by more accurate estimation of the series.

The bounds for the convergence rate can be extended for regenerative processes described in { Section \ref{sec:intro}}.

\paragraph{The work is supported by RFBR grant No~17-01-00633 A.}



 }

\end{document}